\documentclass[12pt,a4paper,leqno]{article}

\usepackage{fullpage,amsmath,amsthm,amsfonts,amscd,graphicx,latexsym,amssymb,eucal,curves}
\usepackage[all]{xy}

\newcommand {\Pfeil} {\longrightarrow}

\newcommand{\abs}[1]{\lvert#1\rvert}

\renewcommand{\phi}{\varphi}

\renewcommand{\epsilon}{\varepsilon}

\newcommand{\cl}{{\rm cl}}

\newcommand{\gen}{{\rm gen}}

\newcommand{\ns}{{\rm ns}}

\newcommand{\sharf}{{\,\sharp\,}}
\newcommand{\sing}{{\rm sing}}

\renewcommand{\top}{{\rm top}}

\DeclareMathOperator{\Char}{char}

\DeclareMathOperator{\PGL}{PGL}
\DeclareMathOperator{\PSL}{PSL}

\newtheorem{Def}{Definition}

\newtheorem{Rem}[Def]{Remark}

\newtheorem{Prop}[Def]{Proposition}

\newtheorem*{Satz*}{Satz}

\newtheorem{Lemma}[Def]{Lemma}

\newtheorem{Ex}{Example}

\newtheorem{Thm}{Theorem}

\newcommand{\A}{\mathbb{A}}

\newcommand{\R}{\mathbb{R}}

\newcommand{\C}{\mathbb{C}}

\newcommand{\F}{\mathbb{F}}

\renewcommand{\P}{\mathbb{P}}

\newcommand {\fU} {\mathfrak U}

\newcommand {\cL} {\mathcal L}

\newcommand {\cX} {\mathcal X}

\begin{document}

\title{On Mumford Orbifolds}
\author{Patrick Erik Bradley}
%
%
\date{}
%
\maketitle
\begin{abstract}
Formulae for the number of branch points of one-dimensional orbifolds
defined over a non-archimedean local field and uniformisable by discrete projective
linear groups are given. They depend only on the uniformising group.
The method of equi\-va\-ri\-ant pasting reveals the possible relative
position of the branch points.

\end{abstract}
%
\section{Introduction}
\label{intro}

Non-archimedean orbifolds first appear in \cite{AndIII}
in the context of uniformising differential equations for certain
covers of the $p$-adic projective line. They are in particular
useful for classifying
triangle groups. Special types of $p$-adic triangle groups,
appearing in characteristic zero only for the residue characteristics
$2$, $3$ and $5$ are Kato's triangle groups {\em of Mumford type} \cite{Kat00}
and \cite{CKK01}.
They belong to one-dimensional orbifolds
 allowing three-point covers of the projective line by
Mumford curves.   

\smallskip
For his uniformisation purpose,
Andr\'e developed a theory of orbifold fundamental groups in characteristic zero which turned out to be very useful in classifying components
of $p$-adic Hurwitz spaces for covers between Mumford curves. This
is the theme of the author's dissertation \cite{Brad02}.

\medskip
Kato's extensive study of triangle groups of Mumford type  yields them all. 
His 
quotient pasting method is  applied in this paper to all finitely
generated discrete projective linear groups over local non-archimedean fields
for finding formulae for the number of branch points of covers defined as
quotients by those groups, and to gain control of their relative position.

\bigskip
The setup for this paper is Berkovich's theory of strictly $K$-analytic
spaces \cite{Ber90}, successfully giving rigid analytic spaces a topology
suitable for studying covers of such spaces. 
For example, Berkovich paths on a one-dimensional manifold $X$ may be 
seen as continuous versions of their rigid analytic counterparts
on the graphs associated to reductions of $X$ along pure affinoid 
coverings of $X$.

\bigskip
Until recently, the author was unaware of an unpublished
article by van der Put and Voskuil \cite{vdPV01} proving a special case of 
these results with entirely different methods.

\bigskip
A large portion of this paper can be found in the author's dissertation
\cite{Brad02} which may serve as an introduction to the theory
of $K$-analytic spaces.

\section{Orbifolds}
\label{Orbifolds}

Let $K$ be a non-archimedean local field. Its residue field
will be denoted $k$.
A {\em manifold} over $K$ or simply a {\em $K$-manifold}, is 
paracompact, strictly $K$-analytic space $X$, admitting
locally an \'etale morphism to $\A_K^{\dim X}$.
A space fulfilling the latter property is  called {\em smooth}.
 There is also an ad-hoc
definition of manifold allowing some non-smooth spaces such as
affinoid discs to be named manifolds \cite{AndIII}.

\begin{Def}
An {\em orbifold} $\cX=(X,(Z_i,G_i))$ over $K$ consists of a 
$K$-manifold $X$ together with a locally finite family $(Z_i)$
of irreducible divisors and finite groups $G_i$ with the property
that for each $i$ there is an open $U_i\supseteq Z_i$ in $X$
and
a ramified Galois cover $V_i\stackrel{/G_i}{\Pfeil} U_i$
un-ramified outside $Z_i$.  
\end{Def}

If $\Char K=0$, this is a special case of an orbifold in the sense
of \cite{AndIII}: the groups $G_i$ may be replaced by cyclic
groups. Andr\'e then replaces these groups by their orders.

\smallskip
If $\Char K>0$, there is a nice subclass of
orbifolds with $G_i$  of the form
$$
G_i\cong B(t,n):=E_t\rtimes C_n,
$$
where $E_t:=C_p^t$ and $n\mid p^t-1$ (depending on $i$, of course).
These
orbifolds are called {\em ordinary}, and may be written as
$\cX=(X,(\xi_i;n_i,t_i)$ with $n_i\mid p^{t_i}-1$.

\medskip
We will always assume that all groups $G_i$ of the orbifold
 $\cX=(X,(Z_i,G_i))$ are non-trivial.

\begin{Def}
A {\em Mumford orbifold} is a one-dimensional orbifold
$\cX=(X,(\xi_i,G_i))$ admitting a Galois cover $Y\to X$
by a Mumford curve $Y$, ramified exactly above the $\xi_i$
with decomposition group isomorphic to $G_i$. 
\end{Def}

The importance of ordinarity is

\begin{Lemma}
In characteristic $p>0$ any Mumford orbifold is ordinary.
\end{Lemma}

\begin{proof}
\cite[Lemma 4.4.2]{Kat02}
\end{proof}

We also have

\begin{Lemma}
If $\cX=(X,(\xi_i,G_i))$ is a Mumford orbifold, then $X$ is a
Mumford curve.
\end{Lemma}

\begin{proof}
For example, the proof of \cite[Lemma 4.2]{Brad02} exhibits this fact.
\end{proof}

\section{Kato trees}
\label{Katotrees}

We will be subsequently working with topological graphs. Only
the edges homeomorphic to the interval $[0,1]$ will be called ''edges''.
The ones homeomorphic to $[0,1)$ are named {\em cusps}.

\medskip
Let $\Omega\subseteq\P^1$ be a connected, open, strictly analytic subset.
It is well-known that $\Omega$ is simply connected. In fact, there
is exactly one injective path connecting two points on $\Omega$.
Thus, for each $x\in\Omega$ we get a partial order $\le_x$ on $\Omega$
by setting
$$
y\le_x z \;:\Leftrightarrow\;\text{$z$ is on the path $x\leadsto y$}.
$$
It is easily seen that for any pair $y,z\in\Omega$ there is a unique
$\sup(y,z)\in\Omega$ for $\le_x$.
This leads to an embedding of $\Omega$ into the space $\hat\Omega$
of all partial orders with suprema for pairs compatible with paths
on $\Omega$ \cite[4.1.3]{Ber90}.

\begin{Def}
The connected hull of $\hat\Omega\setminus\Omega$ in $\hat\Omega$ 
is called the 
{\em skeleton} of $\Omega$.
\end{Def}

Let $X$ be a non-singular projective curve over $K$. Replacing $K$
by a finite extension, if necessary, guarantees the existence of
a pure affinoid covering $\fU$ of $X$ such that the reduction
$\overline{X_{\fU}}$ along $\fU$ is semi-stable \cite{vdP84}.
Let $\pi_{\fU}\colon X\to\overline{X_{\fU}},\;x\mapsto \bar{x}$ 
be the corresponding
Tate map.
It is well known that 
$$
\pi_{\mathfrak U}^{-1}(\bar{x})\text{ is a(n) }
\begin{cases}
	\text{point}, &\text{if }\bar{x}\in \bar{X}_\gen\\
	\text{open disk}, &\text{if }\bar{x}\in\bar{X}_\ns \\
	\text{open annulus}, &\text{if }\bar{x}\in\bar{X}_\sing 
\end{cases}					
$$
\cite[4.3.1]{Ber90}.
Here $\bar{X}_\gen$, $\bar{X}_\ns$ resp.~$\bar{X}_\sing$ means the
set of generic, non-singular resp.~singular points of
the connected curve $\bar{X}$ defined over $k$.
The skeleton $\Delta(\pi_{\fU}^{-1}(\bar{x}))$ is therefore either empty
($\bar{x}\in\bar{X}_\gen$ or $\bar{X}_\ns$) or homeomorphic to the open
unit interval.
The graph $\Delta_{\fU}(X)$ is defined as follows:
its vertex set is $\pi_{\fU}^{-1}(\bar{X}_\gen)$, and the
edges are the closures $\pi_{\fU}^{-1}(\bar{x})^\cl$ in $X$
(and therefore homeomorphic either to $[0,1]$ or $[0,1)$). 

\begin{Def}
The graph $\Delta_{\fU}(X)$ is called the {\em analytic skeleton} of $X$ with
respect to the pure affinoid covering $\fU$.
\end{Def}

\smallskip 
A discrete group
$G$ acting discontinuously on $\Omega$ induces an action on
the skeleton $\Delta(\Omega)$, which in fact is a tree \cite[4.1.7]{Ber90}.
If $X=\Omega/G$ is a non-singular irreducible projective curve, then there
is a pure affinoid covering $\fU$ of $X$ such that
$$
\Delta(\Omega)/G\cong\Delta_{\fU}(X).
$$
The topological graph of groups $(\Delta_{\fU}(X),G_\bullet)$ will be called
a {\em quotient skeleton}.

\smallskip
Let now be given a finitely generated discrete group $N\subseteq\PGL_2(K)$.
If we take $\Omega^*(N)$ to be the complement of the set 
$\cL^*(N)$
of limit points
of $N$ and the 
fixed points of $N$'s elements of finite order, then its skeleton
$\Delta^*(N):=\Delta(\Omega^*)$ is a tree whose set of  cusps 
is exactly $\cL^*(N)$.

\begin{Def}
The quotient skeleton $\Gamma^*(N):=(\Delta^*(N)/N,N_\bullet)$
is called the {\em Kato graph} of $N$.
\end{Def} 

There is also $\Omega(N):=\P^1\setminus\cL(N)$,
the complement of the set $\cL(N)$ of limit points of $N$ in $\P^1$. The
corresponding quotient skeleton 
$$
\Gamma(N)=(\Delta(\Omega(N))/N,N_\bullet) 
$$
is a finite graph of groups, and the map $\Omega(N)\to\Omega(N)/N$
is an orbifold cover of the Mumford curve $\Omega(N)/N$. 
The graph $\Gamma(N)$ is obtained by a 
{\em contraction} of the Kato graph $\Gamma^*(N)$:
cusps are cut off, and edges $e$ whose stabiliser equals either
$N_{o(e)}$ or $N_{t(e)}$ are replaced by the extremal vertex $v$
with the larger group, if the valency of $v$ is less than $3$.

\smallskip
Kato states in \cite[Proposition 2.2]{Kat00}

\begin{Lemma}
There is a one-to-one correspondence between the set of cusps
of $\Gamma^*(N)$ and the set of branch points of the cover
$\Omega(N)\to\Omega(N)/N$. The stabiliser of a cusp is
the decomposition group of the corresponding branch point.
\end{Lemma}

The relationship between Mumford orbifolds and Kato graphs is

\begin{Thm}
A one-dimensional orbifold is a Mumford orbifold if and only if 
there exists a global Galois orbifold chart whose quotient skeleton
is a Kato graph.
\end{Thm}

\begin{proof}
Denote the orbifold with $\cX=(X,(\xi,e_i))$.

\smallskip
The implication $\Rightarrow$ is quite clear: for an orbifold chart 
$\phi\colon C\stackrel{/G}{\Pfeil}\cX$ 
take $\omega\colon\Omega\to C$ to be the universal cover of the Mumford curve $C$.
Then the composition $\phi\circ\omega$ is the quotient by an action
of a discrete group $N\subseteq\PGL_2(K)$, and the quotient graph
$(\Delta^*(N)/N,N_\bullet)$
is a Kato graph.

\smallskip
For $\Leftarrow$: let $C\to\cX$ be a chart whose quotient skeleton is a Kato graph for $N$. Being a fundamental group of a graph whose vertex and
edge stabilisers are finite subgroups of $\PGL_2(K)$,
$N$ contains a free normal subgroup $H$ of finite index \cite[I.\S 3]{GvP80}.
Let $\Omega\subseteq\P^1$ be their domain of regularity. Then there is a 
commuting diagramme
$$
\xymatrix{
\Omega\ar_\top^{/H}[r] \ar_{/N}[dr] & Y \ar[d]\\
 &\cX
}
$$
with a Mumford curve $Y$. Since the vertical arrow is an orbifold chart,
we are done.
\end{proof}

In his construction of triangle groups, Kato 
 first studies the Kato trees
of the finite subgroups of $\PGL_2(K)$ \cite{Kat00}.

\begin{Def}
The Kato trees $T^*(G)$ for the finite groups of $\PGL_2(K)$ are
called {\em elementary Kato trees}. 
\end{Def}

\begin{Prop}
The elementary Kato trees look like this, if $\Char K=0$:

\begin{enumerate}
\item If $\Char k>5$, then all elementary Kato trees have only one vertex:

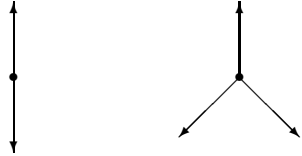
\setlength{\unitlength}{.5cm}
\begin{picture}(15,6)
\put(6,3){\circle*{.2}}
\put(6,3){\vector(0,1){2}}
\put(6,3){\vector(0,-1){2}}

\put(12,3){\circle*{.2}}
\put(12,3){\vector(0,1){2}}
\put(12,3){\vector(1,-1){1.6}}
\put(12,3){\vector(-1,-1){1.6}}
\end{picture}

\item
If $\Char k\le 5$, then the elementary Kato trees for the
groups of order divisible by $p$ are also star-shaped as in the above case, but
many of them
are instable (after contracting the cusps). 
\end{enumerate}
\end{Prop}

\begin{proof}
\cite{Kat00}.
\end{proof}

\begin{Rem}
The existence of edges 
in $\Char k\le 5$ turns out to be 
important for constructing $p$-adic triangle groups of Mumford type
{\rm\cite{Kat00}}.
\end{Rem}

If $\Char K>0$, the Kato tree for $G=C_p^t$, as defined above, is empty
\cite[4.1.4.(i)]{Ber90}.
We will resolve the problem in the following manner:
since all $\gamma\in C_p^t\setminus\{1\}$ are parabolic, the
fixed point set in $\P^1$ of $G$ is  a one-cusped subtree of 
$\Delta(\P^1\setminus\P^1(K))$,
and may be contracted to a vertex with one cusp.
We will call that tree, by abuse of language, the Kato tree for $C_p^t$. 

\medskip
The geodesic in $\P^1$ between the two $K$-rational fixed points 
of an elliptic transformation $\gamma$
is called the {\em mirror} $M(\gamma)$ of $\gamma$. 

\medskip
Let $P(2,p^t)$ denote either $\PGL_2(\F_{p^t})$ or $\PSL_2(\F_{p^t})$.
In the former case, we will abbreviate
$n_-:=p^t-1$ and $n_+:=p^t+1$, whereas in the latter,
$n_-:=\frac{p^t-1}{2}$ and $n_=:=\frac{p^t+1}{2}$.

\medskip
In the proof of the following proposition, we will take particular care
whether (two-pointed) mirrors are  {\em folded}, i.e.~mapped onto a 
half-line, or not.

\begin{Prop}	\label{elementary p}
If $\Char K=p>0$, then any elementary Kato tree $T^*(G)$
for $G\subseteq\PGL_2(K)$ viewed as a subgroup of $\PGL_2(\F_{p^m})$ is one of
the following:
\begin{enumerate}

	\item $G=C_n$ for $(n,p)=1\colon$
\setlength{\unitlength}{.5cm}
\begin{picture}(6,0)
\put(1,-.2){$n$}
\put(4,0){\vector(-1,0){2}}
\put(4,0){\vector(1,0){2}}
\put(4,0){\circle*{.2}}
\put(6.5,-.2){$n$}
\put(3.6,.3){$C_n$}
\end{picture}
	\item $G=D_n$ for 
\begin{enumerate}
	\item $p\neq 2$ and $n\mid p^m\pm 1\colon$
\setlength{\unitlength}{.5cm}
\begin{picture}(8,5)
\put(4,2){\circle*{.2}}
\put(4,2){\vector(-1,-1){1.5}}
\put(4,2){\vector(1,-1){1.5}}
\put(4,2){\vector(0,1){2}}
\put(4.3,1.8){$D_n$}
\put(2,0){$2$}
\put(5.7,0){$2$}
\put(3.8,4.2){$n$}
\end{picture}

	\item $p=2$ and $(n,2)=1\colon$
\setlength{\unitlength}{.5cm}
\begin{picture}(6,2)
\put(1,-.2){$2$}
\put(4,0){\vector(-1,0){2}}
\put(4,0){\vector(1,0){2}}
\put(4,0){\circle*{.2}}
\put(6.5,-.2){$n$}
\put(3.6,.3){$D_n$}
\end{picture}

\end{enumerate} 

	\item $G=B(t,n)$ for $t\le m$ and
\begin{enumerate}
	\item  $n\mid p^m-1$, $n\mid p^t-1$, $n>1\colon$
\setlength{\unitlength}{.5cm}
\begin{picture}(6,2)
\put(1,-.2){$n$}
\put(4,0){\vector(-1,0){2}}
\put(4,0){\vector(1,0){2}}
\put(4,0){\circle*{.2}}
\put(6.5,-.2){$B(t,n)$}
\put(3.1,.3){$B(t,n)$}
\end{picture}
	\item $n=1\colon$
\setlength{\unitlength}{.5cm}
\begin{picture}(15,2)
\put(4,0){\circle*{.2}}
\put(4,0){\vector(-1,0){2}}
\put(1,-.2){$E_t$}
\put(4.2,-.2){$E_t$}
\end{picture}

\end{enumerate}
	\item $G=P(2,p^t)\colon$
\setlength{\unitlength}{.5cm}
\begin{picture}(6,2)
\put(1,-.2){$n_+$}
\put(4,0){\vector(-1,0){2}}
\put(4,0){\vector(1,0){2}}
\put(4,0){\circle*{.2}}
\put(6.5,-.2){$B(t,n_-)$}
\put(2.8,.3){$P(2,p^t)$}
\end{picture}

	\item $G=T$ for $p\neq 2,3\colon$
\setlength{\unitlength}{.5cm}
\begin{picture}(8,5)
\put(4,2){\circle*{.2}}
\put(4,2){\vector(-1,-1){1.5}}
\put(4,2){\vector(1,-1){1.5}}
\put(4,2){\vector(0,1){2}}
\put(4.3,1.8){$T$}
\put(2,0){$3$}
\put(5.7,0){$3$}
\put(3.8,4.2){$2$}
\end{picture}

	\item $G=O$ for $p\neq 2,3\colon$
\setlength{\unitlength}{.5cm}
\begin{picture}(8,5)
\put(4,2){\circle*{.2}}
\put(4,2){\vector(-1,-1){1.5}}
\put(4,2){\vector(1,-1){1.5}}
\put(4,2){\vector(0,1){2}}
\put(4.3,1.8){$O$}
\put(2,0){$2$}
\put(5.7,0){$3$}
\put(3.8,4.2){$4$}
\end{picture}

	\item $G=I$ for $5\mid p^{2m}-1$ and
\begin{enumerate}
	\item  $p\neq 2,3,5\colon$
\setlength{\unitlength}{.5cm}
\begin{picture}(8,5)
\put(4,2){\circle*{.2}}
\put(4,2){\vector(-1,-1){1.5}}
\put(4,2){\vector(1,-1){1.5}}
\put(4,2){\vector(0,1){2}}
\put(4.3,1.8){$I$}
\put(2,0){$2$}
\put(5.7,0){$3$}
\put(3.8,4.2){$5$}
\end{picture}

	\item $p=3\colon$
\setlength{\unitlength}{.5cm}
\begin{picture}(6,2)
\put(1,-.2){$5$}
\put(4,0){\vector(-1,0){2}}
\put(4,0){\vector(1,0){2}}
\put(4,0){\circle*{.2}}
\put(6.5,-.2){$B(1,2)$}
\put(3.8,.3){$I$}
\end{picture}

\end{enumerate}
\end{enumerate}
\end{Prop}

\begin{proof}
The list of groups is Dickson's classification of finite subgroups of
$\PGL_2(K)$ \cite[II.8.27]{Hu79}, to be found also in \cite{VM80},
\cite[Theorem 2.9]{CKK01}.
The trees of 1., 2.(a), 5., 6., 7.(a) have been constructed in
\cite[Appendix A]{Kat00} for $\Char K=0$. Kato's construction is also
possible in our
situation. 
3.(b) is our convention from above.

\smallskip
2.(b): $D_n=C_n\rtimes C_2$, and $C_2$ interchanges the two $K$-rational
fixed points of $C_n$, thus folding the mirror of a generator of
$C_n$, while
the generator of $C_2$ is parabolic. 

\smallskip
3.(a): The tree is the intersection of the two mirrors $M(\sigma)$ 
and $M(\gamma)$ of an elliptic $\sigma$ of order $n>1$ and $K$-rational
fixed points, say, $0$ and $\infty$, and a parabolic $\gamma$
with fixed point $0$. So one of the cusps of $M(\sigma)$ is stabilised
by all of $B(t,n)$.

\smallskip
4.: $P(2,q^t)$ contains the dihedral group $D_{n_+}$ \cite{VM80}.
The element of order two in it folds the mirror of any other generator
of order $n_+$. 
Since $T^*(P(2,q^t)$ has only two cusps, anyway, the other one
is stabilised by $B(t,n_-)$. Since, according to \cite[Lemma 4.3]{CKK01},
every path emanating from the vertex $v$ with stabiliser $P(2,p^t)$
is locally stabilised by either $C_{n_+}$ or $B(t,n_-)$, 
the Kato tree has the shape
as drawn.   

\smallskip
7.(b): As $I$ contains dihedral groups $D_5$ and $D_2$, the mirrors of elements
of orders $2$ and $5$ are folded. But the former mirror has a cusp
in common with the mirror of a parabolic transformation.
\end{proof}

\begin{Def}
A tree of groups with non-trivial edge groups is called an {\em irreducible
tree}. A maximal irreducible subtree of a graph $\Gamma$ of groups is called
an {\em irreducible component} of $\Gamma$.
\end{Def}

The following Proposition is decisive for the rest of this paper.

\begin{Prop}	\label{glue}
Every irreducible Kato tree $T=(T^*(N),N_\bullet)$ is obtained by 
glueing the elementary Kato trees $T^*(N_v)$ for the vertex groups
of $T$ along the trees $T^*(N_e)$ for the edge groups:
$$
T\cong\lim\limits_{\longrightarrow}T^*(N_\bullet).
$$
\end{Prop}

\begin{proof}
The proof in \cite[Proposition 3.9]{Kat00} works in all characteristics.
\end{proof}

\section{The structure of Kato graphs}
\label{structure}

An important question in the analysis
of Kato graphs is whether mirrors  fold or not under the action of
a given group. 

\subsection{The characteristic zero case}
\label{zero}

\begin{Thm}	\label{formula}
Let $\Char K=0$. If $C$ denotes the number of cyclic vertex groups,
$c$ the number of cyclic edge groups, $D$ the number of non cyclic
vertex groups, and $d$ the number of non cyclic edge groups of $\Gamma(N)$,
then the number of cusps of $\Gamma^*(N)$ is given by
$$
\#\partial\Gamma^*(N)=3(D-d)+2(C-c).
$$
\end{Thm}

\begin{proof}
First let $\Gamma^*(N)=T$ be an irreducible Kato tree.
Then, according to Proposition \ref{glue},
$T\cong\lim\limits_{\longrightarrow}T^*(N_\bullet)$ 
is obtained
by glueing the elementary Kato trees for the vertex groups
along the trees for the edge groups.
So, all we have to do is  to examine the glueing process, and check the
formula at each step. A detailed description of the glueing morphisms
as quotients of the injective glueing morphisms 
$\Delta^*(N_e)\to\Delta^*(N_v)$ for $v$ an extremity of $e$
has already been done in \cite{Kat00}.

\smallskip
If $\Char k>5$, then only cyclic groups can occur as edge groups.
\cite{Kat00} tells us that in this case segments are obtained somewhat
like

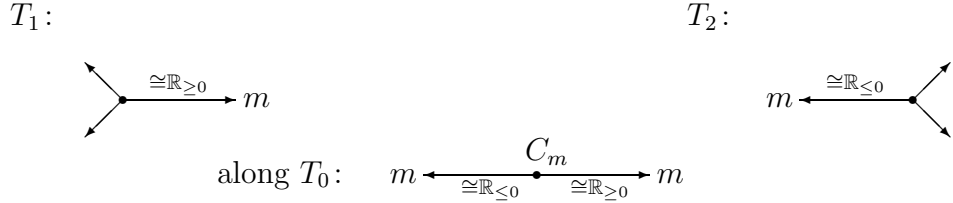
\setlength{\unitlength}{.5cm}
\begin{picture}(26,8)
\put(1,7){$T_1\colon$}
\put(4,5){\circle*{.2}}
\put(4,5){\vector(-1,1){1}}
\put(4,5){\vector(-1,-1){1}}
\put(4,5){\vector(1,0){3}}
\put(7.2,4.8){$m$}
\put(4.7,5.3){$\scriptstyle\cong\R_{\ge 0}$}

\put(19,7){$T_2\colon$}
\put(25,5){\circle*{.2}}
\put(25,5){\vector(-1,0){3}}
\put(25,5){\vector(1,1){1}}
\put(25,5){\vector(1,-1){1}}
\put(21.1,4.8){$m$}
\put(22.7,5.3){$\scriptstyle\cong\R_{\le 0}$}

\put(6.5,2.8){\text{along $T_0\colon$}}
\put(15,3){\circle*{.2}}
\put(15,3){\vector(-1,0){3}}
\put(15,3){\vector(1,0){3}}
\put(11.1,2.8){$m$}
\put(18.2,2.8){$m$}
\put(14.7,3.4){$C_m$}
\put(13,2.5){$\scriptstyle\cong\R_{\le 0}$}
\put(15.9,2.5){$\scriptstyle\cong\R_{\ge 0}$}
\end{picture} 

with glueing morphisms
$$
\phi_1\colon T_0\to T_1, \; t\mapsto \abs{t},
\qquad \phi_2\colon T_0\to T_1,\;t\mapsto -\abs{1-t},
$$
resulting in

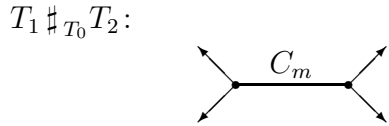
\begin{picture}(4,6)
\put(1,4.5){$T_1\sharf_{T_0}T_2\colon$}
\put(7,3){\circle*{.2}}
\put(10,3){\circle*{.2}}
\put(7,3){\line(1,0){3}}
\put(7,3){\vector(-1,1){1}}
\put(7,3){\vector(-1,-1){1}}
\put(10,3){\vector(1,1){1}}
\put(10,3){\vector(1,-1){1}}
\put(7.9,3.3){$C_m$}
\end{picture}

This proves the formula for irreducible Kato trees, if $\Char k>5$.

\smallskip
Now, for $\Char k\le 5$, there are elementary Kato trees with
edges (homeomorphic to $[0,1]$). 
This makes it possible to glue along trees for non cyclic
groups in order to obtain trees with edges. Kato has shown \cite{Kat00}
that in those cases one of the following occurs:

\begin{itemize}
\item Both glueing morphisms are injective.
\item One of the glueing morphisms folds an unfolded mirror of $T_0$.
\item Both glueing morphisms fold an unfolded mirror of $T_0$.
\end{itemize}
Any of these possibilities reduces the total amount of ends by three,
and the formula holds in this case.

\medskip
Now that the formula has been proven for irreducible trees, we examine
the case that $T=\Gamma^*(N)$ is a general Kato tree. 
Then $N$ is a free tree product of fundamental groups $N_i$ of 
irreducible Kato trees $T_i=\Gamma^*(N_i)$ which are connected
within $T$ by edges with trivial stabiliser. Obviously,
the number of cusps of $T$ equals the sum of the number of cusps belonging to
the irreducible 
 components $T_i$. The formula is therefore valid for all Kato trees.

\medskip
Now, let $\Gamma=\Gamma^*(N)$ be a Kato graph of positive genus $g$.
Since the fundamental group $N$ of $\Gamma$ contains a free normal subgroup
of rank $g$, we
 may remove $g$ edges with trivial stabiliser in such a way as to obtain a tree
of groups with the same number of vertices.
This proves the
 formula in all genera of $\Gamma^*(N)$.
\end{proof}

\begin{Ex}
One of the simpler triangle groups is for $\Char k=5$:
$T_0=\Gamma^*(D_5)$, $T_1=\Gamma^*(A_5)$, $T_2=\Gamma^*(D_{10m})$.

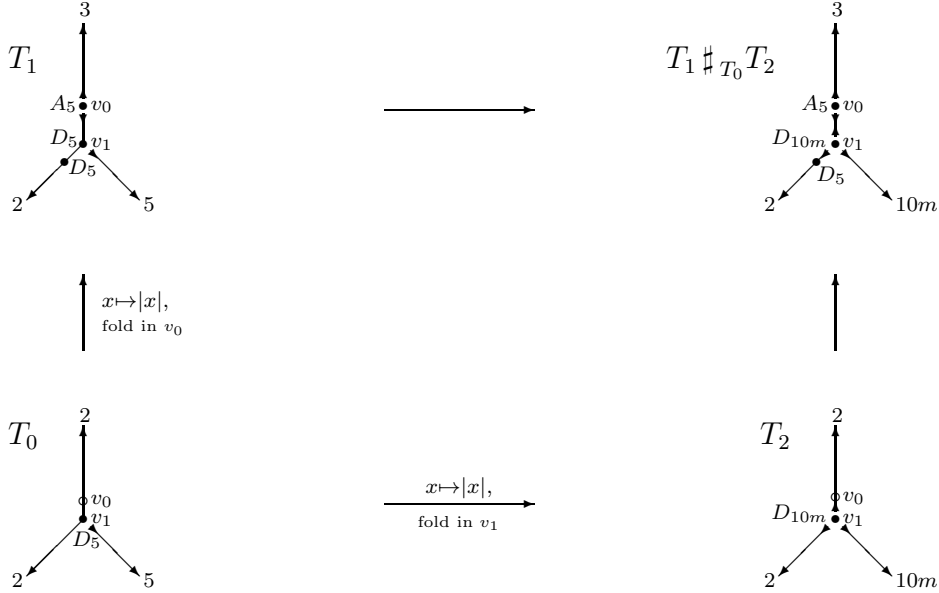
\setlength{\unitlength}{.5cm}
\begin{picture}(20,17)

\put(3,2.5){\vector(0,1){2.5}}
\put(2.9,5.1){$\scriptstyle 2$}
\put(3,2.5){\circle*{.2}}
\put(3,3){\circle{.2}}
\put(3,2.5){\vector(-1,-1){1.5}}
\put(3.2,2,3){\vector(1,-1){1.3}}
\put(3.2,2.3){\vector(1,-1){.2}}
\put(1.1,.7){$\scriptstyle 2$}
\put(4.6,.7){$\scriptstyle 5$}
\put(2.7,1.8){$\scriptstyle D_5$}
\put(3.2,2.4){$\scriptstyle v_1$}
\put(3.2,2.9){$\scriptstyle v_0$}
\put(1,4.5){$T_0$}

\put(3,7){\vector(0,1){2}}
\put(3.5,7.5){\shortstack[l]{$\scriptstyle x\mapsto\abs{x},$\\
			\text{\rm \tiny fold in $v_0$}}}

\put(3,13.7){\vector(0,1){2}}
\put(3,13.8){\vector(0,1){.2}}

\put(3,13.3){\line(0,-1){.8}}
\put(3,13.2){\vector(0,-1){.2}}

\put(2.9,15.9){$\scriptstyle 3$}
\put(3,12.5){\circle*{.2}}
\put(3,13.5){\circle*{.2}}
\put(3,12.5){\vector(-1,-1){1.5}}
\put(3.2,12,3){\vector(1,-1){1.3}}
\put(3.2,12.3){\vector(1,-1){.2}}
\put(1.1,10.7){$\scriptstyle 2$}
\put(4.6,10.7){$\scriptstyle 5$}
\put(2.6,11.7){$\scriptstyle D_5$}
\put(3.2,12.4){$\scriptstyle v_1$}
\put(3.2,13.4){$\scriptstyle v_0$}
\put(1,14.5){$T_1$}
\put(2.1,13.4){$\scriptstyle A_5$}
\put(2.1,12.5){$\scriptstyle D_5$}
\put(2.5,12){\circle*{.2}}

\put(23,2.7){\vector(0,1){2.3}}
\put(23,2.8){\vector(0,1){.2}}
\put(22.9,5.1){$\scriptstyle 2$}
\put(23,2.5){\circle*{.2}}
\put(23,3.1){\circle{.2}}
\put(22.8,2.3){\vector(-1,-1){1.3}}
\put(22.8,2.3){\vector(-1,-1){.2}}
\put(23.2,2,3){\vector(1,-1){1.3}}
\put(23.2,2.3){\vector(1,-1){.2}}
\put(21.1,.7){$\scriptstyle 2$}
\put(24.6,.7){$\scriptstyle 10m$}
\put(21.3,2.5){$\scriptstyle D_{10m}$}
\put(23.2,2.4){$\scriptstyle v_1$}
\put(23.2,3){$\scriptstyle v_0$}
\put(21,4.5){$T_2$}

\put(23,7){\vector(0,1){2}}

\put(23,13.7){\vector(0,1){2}}
\put(23,13.8){\vector(0,1){.2}}

\put(23,13.3){\line(0,-1){.6}}
\put(23,13.2){\vector(0,-1){.2}}
\put(23,12.8){\vector(0,1){.2}}

\put(22.9,15.9){$\scriptstyle 3$}
\put(23,12.5){\circle*{.2}}
\put(23,13.5){\circle*{.2}}
\put(22.8,12.3){\vector(-1,-1){1.3}}
\put(22.8,12.3){\vector(-1,-1){.2}}
\put(23.2,12,3){\vector(1,-1){1.3}}
\put(23.2,12.3){\vector(1,-1){.2}}
\put(21.1,10.7){$\scriptstyle 2$}
\put(24.6,10.7){$\scriptstyle 10m$}
\put(22.5,11.5){$\scriptstyle D_5$}
\put(23.2,12.4){$\scriptstyle v_1$}
\put(23.2,13.4){$\scriptstyle v_0$}
\put(18.5,14.5){$T_1\sharf_{T_0}T_2$}
\put(22.1,13.4){$\scriptstyle A_5$}
\put(21.3,12.5){$\scriptstyle D_{10m}$}
\put(22.5,12){\circle*{.2}}

\put(11,13.4){\vector(1,0){4}}
\put(11,2.9){\vector(1,0){4}}
\put(12.1,3.2){$\scriptstyle x\mapsto \abs{x},$}
\put(11.9,2.3){\text{\rm \tiny fold in $v_1$}}
\end{picture}

Note that the point $v_0$ in $T_1$ is a vertex, whereas in
$T_0$ and $T_2$ it is not: it is a marked point $\circ$ on the cusp
viewed as the half-open interval $[0,1)$.
\end{Ex}

\subsection{The general case}
\label{general}

Let now $K$ be of any characteristic. 
If $\Char K>0$, then the formula of Theorem \ref{formula}
has to be stated in more general terms, as Kato trees for non-cyclic groups
may have different numbers of cusps. Let $V$ denote the set of vertices
and $E$ the set of edges of  a given graph.

\begin{Thm}
The number of cusps in a Kato graph $\Gamma^*(N)$ is given by
$$
\#\partial\Gamma^*(N)=\sum\limits_{v\in V}\#\partial T^*(N_v)
-\sum\limits_{e\in E}\#\partial T^*(N_e).
$$
\end{Thm}

\begin{proof}
The proof reduces to considering $\Char K>0$.
As in the characteristic zero case, we need to treat irreducible
Kato trees, only.
In positive characteristic, the glueing procedure turns out to be 
of the following:
\begin{enumerate}
\item For edge groups $N_e$ we have $\#\partial T^*(N_e)\le 2$.
\item If $\#\partial T^*(N_e)=1$, then both glueing morphisms
$T^*(N_e)\to T^*(N_{o(e)})$ and $T^*(N_e)\to T^*(N_{t(e)})$
are injective.
\item If $\#\partial T^*(N_e)=2$, then the glueing morphisms
either both fold the line $ T^*(N_e)$ in two different points,
or one of the glueing morphisms is an isomorphism of trees
with injective  morphisms between stabilisers. 
\end{enumerate}
The proof of these claims.

\smallskip
1. Edge stabilisers are Borel groups \cite[Lemma 4.1]{CKK01}
(The cited Lemma holds without the restrictions of  \cite[Section 4]{CKK01}). 

\smallskip
2. Obvious.

\smallskip
3. We have $N_e=B(t,n)$ with $n>1$.
If $N_e$ is cyclic (of order necessarily prime to $p$), then glueing
is as in the characteristic zero case.

\smallskip
Let now $N_e=B(t,n)$ with $t>0$. There are two possibilities:
either the stabiliser of the origin $o(e)$ is of Borel type or not.
In the first case, $N_{o(e)}=B(t',n')$,
then $n'=n$ because the prime-to-$p$ part of $N_e$ is a maximally
cyclic subgroup of $N_{o(e)}$ \cite[Lemma 1]{Her79}. 

\smallskip
In the second case, we have $N_{o(e)}=P(2,p^{t'})$. Then, as 
we have seen in the proof of Proposition 
\ref{elementary p}.~4., $t'=t$ and $n=n_-$.
An element of order two in $P(2,p^t)$ folds the geodesic tree $T^*(B(t,n_-)$
in its vertex $v_0$.
The glueing morphism maps $v_0$ somewhere onto the cusp
of $T^*(N_{o(e)})$ with stabiliser $B(t,n)$. The image $v_0'$ of $v_0$ in turn
will become $t(e)$ in $T^*(N)$. 
The pre-image of $e$ in $T^*(N_{o(e)})$ is the path from its vertex
to $v_0'$, and its stabiliser is $B(t,n_-)$.
As shown in the proof of \cite[4.6]{CKK01}, the ends of any
 geodesic in $T^*(N)$ going through $e$ are stabilised by $C_{n_+}$ and
by $B(td,n_-)$, respectively, and
the stabilisers of  points on the latter end form an increasing sequence of
Borel groups, as one moves away from $t(e)$. 
Therefore, $N_{t(e)}=B(t',n_-)$
with $t\mid t'$, and the map $T^*(N_e)\to T^*(N_{t(e)})$ merely replaces
$B(t,n_-)$ by $B(t',n_-)$.
 
\smallskip
These are all possibilities.
\end{proof}

\begin{Ex} A non-elementary Kato tree with two cusps may be constructed
in this way: let $n=p^t-1$ and $t\mid s$, then

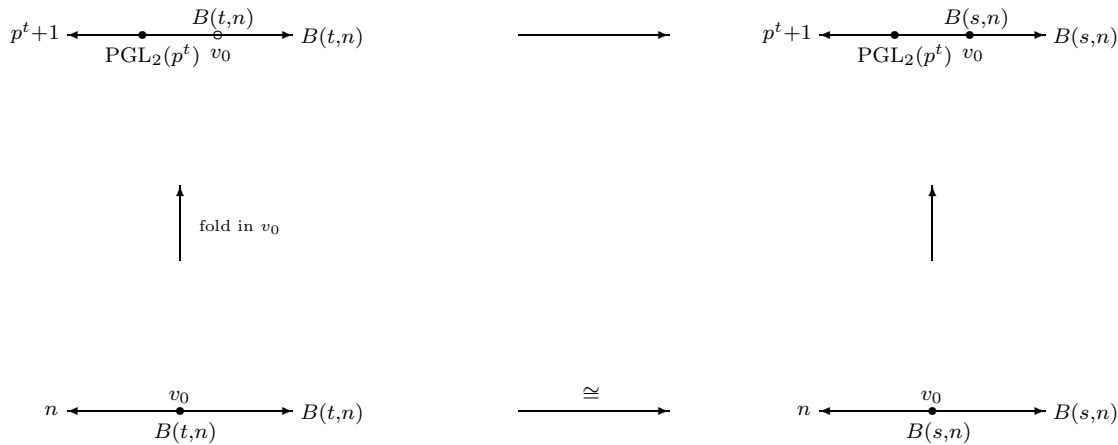
\setlength{\unitlength}{.5cm}

\begin{picture}(20,14)
\put(6,2){\circle*{.2}}
\put(5,2){\vector(-1,0){2}}
\put(5,2){\vector(1,0){4}}
\put(2.4,1.9){$\scriptstyle n$}
\put(9.2,1.8){$\scriptstyle B(t,n)$}
\put(5.7,2.3){$\scriptstyle v_0$}
\put(5.3,1.3){$\scriptstyle B(t,n)$}

\put(26,2){\circle*{.2}}
\put(25,2){\vector(-1,0){2}}
\put(25,2){\vector(1,0){4}}
\put(22.4,1.9){$\scriptstyle n$}
\put(29.2,1.8){$\scriptstyle B(s,n)$}
\put(25.7,2.3){$\scriptstyle v_0$}
\put(25.3,1.3){$\scriptstyle B(s,n)$}

\put(5,12){\circle*{.2}}
\put(5,12){\vector(-1,0){2}}
\put(5,12){\vector(1,0){4}}
\put(7,12){\circle{.2}}
\put(1.5,11.9){$\scriptstyle p^t+1$}
\put(9.2,11.8){$\scriptstyle B(t,n)$}
\put(6.8,11.4){$\scriptstyle v_0$}
\put(4,11.3){$\scriptstyle \PGL_2(p^t)$}
\put(6.3,12.3){$\scriptstyle B(t,n)$}

\put(25,12){\circle*{.2}}
\put(25,12){\vector(-1,0){2}}
\put(25,12){\vector(1,0){4}}
\put(27,12){\circle*{.2}}
\put(21.5,11.9){$\scriptstyle p^t+1$}
\put(29.2,11.8){$\scriptstyle B(s,n)$}
\put(26.8,11.4){$\scriptstyle v_0$}
\put(26.3,12.3){$\scriptstyle B(s,n)$}
\put(24,11.3){$\scriptstyle \PGL_2(p^t)$}

\put(15,2){\vector(1,0){4}}
\put(15,12){\vector(1,0){4}}
\put(6,6){\vector(0,1){2}}
\put(26,6){\vector(0,1){2}}

\put(16.7,2.3){$\scriptstyle \cong$}
\put(6.5,6.8){\text{\rm \tiny fold in $v_0$}}

\end{picture} 

is Cartesian. The stabiliser of the edge in the upper right tree
is the intersection of the two vertex groups:
$\PGL_2(p^t)\cap B(x,n)= B(t,n)$.
\end{Ex}

From the proofs of the two preceding theorems we get

\begin{Prop}
Let $\Gamma$ be a Kato graph. Then
\begin{enumerate}
\item
Each vertex has at most three cusps or edges with non trivial stabilisers 
going out of it.
\item
The vertex stabilisers are generated by the outgoing edges' or cusps'
 groups.
\end{enumerate}
\end{Prop}

\subsection{Separating the branch points}
\label{position}

The knowledge we have over the glueing process allows us
to control the separation of branch points.

\begin{Thm}
If $\mathcal{X}=(X,(x_i,G_i))$ is a Mumford orbifold, then
there is a pure affinoid covering of $X$ separating the points
$x_i$ into triplets, pairs and singlets. 
\end{Thm}

\begin{proof}
An orbifold chart 
$\phi\colon Y\stackrel{/G}{\longrightarrow} \mathcal{X}$ 
with a Mumford curve $Y$ gives us a Kato graph
$\Gamma^*(N)=(\Gamma,G_\bullet)$ with fundamental group 
$N=\lim\limits_{\longrightarrow} G_\bullet$.
As we have seen in the proof of Theorem \ref{formula},
$\Gamma$ is obtained by glueing elementary Kato trees.
Such trees have at most three cusps, and therefore,
at most three cusps emanate from each vertex of $\Gamma$.

\smallskip
Now, the cusps themselves are Berkovich paths from the 
branch points $x_i$ to the nearest vertex, i.e.~a generic point
of a disc with radius in $\abs{K^\times}$ containing all points $x_i$ 
 to which the  cusps emanating from that vertex lead.
In the case of two or three cusps, that disc is the smallest
one containing the corresponding points which are equidistantly
positioned on  the disc's boundary.
Any nearest  vertex in $\Gamma$ with cusps is a 
generic point of a disc with radius in 
$\abs{K^\times}$ containing the corresponding points $x'_i$ in its
boundary. The shortest path to the first vertex gives the distance
between those two generic points. All edges emanating from those
two vertices define two affinoid subsets of $X$ separating the
set  $\{x_i\}$ from $\{x'_i\}$ and the rest
of the marked points of $\mathcal{X}$. As the intersection of the
two affinoid sets consists of boundary components (coming from
edges connecting the two points), we see that it is a pure subset
of both affinoid pieces. Thus, one obtains a pure
affinoid covering fulfilling the requirements of the theorem. 
\end{proof}

%
%
%
%

\begin{thebibliography}{X}

\bibitem{AndIII}Y.~ Andr\'e. {\em $p$-adic orbifolds and monodromy}, chapter III from {\em Period mappings and differential equations. From $\C$ to $\C_p$}, T\^{o}hoku-Hokkaid\^{o} Lectures in Arithmetic Geometry, preprint {\tt math.AG/0203194}

\bibitem{Ber90}V.G.~Berkovich. {\em Spectral Theory and Analytic Geometry over Non-Archimedean Fields}, Mathematical Surveys and Monographs, Number 33, AMS (1990)



\bibitem{Brad02}P.E.~Bradley. {$p$-adische Hurwitzr\"aume}, dissertation, Universit\"at Karlsruhe (2002)


\bibitem{CKK01}G.~Cornelissen, F.~Kato, A.~Kontogeorgis. {\em Discontinuous groups in positive characteristic and automorphisms of Mumford curves}, Math.\ Ann.\ 320, 55-85 (2001)







\bibitem{GvP80}L.~Gerritzen, M.~van der Put. {\em Schottky Groups and Mumford Curves}, Lecture Notes in Mathematics 817, Springer-Verlag, 1980


\bibitem{Her79}F.~Herrlich. {\em Die Ordnung der Automorphismengruppe einer $p$-adischen Schottkykurve}, Math.\ Ann.\ 246, Nr.\ 2, 126-130 (1979/80)





\bibitem{Hu79}B.~Huppert. {\em Endliche Gruppen I}, Grundl.\ d.\ math.\ Wiss.\ 134, Springer-Verlag, 1979

\bibitem{Kat00}F.~Kato. {\em $p$-adic Schwarzian triangle groups of Mumford type}, preprint {\tt math.AG/9908174}

\bibitem{Kat02}F.~Kato. {\em Non-archimedean orbifolds covered by Mumford curves}, preprint (2001) 















\bibitem{VM80}R.~C.~Valentini, M.~L.~Madan. {\em A Hauptsatz of L.\ E.\ Dickson and Artin-Schreier extensions}, J.\ Reine Angew.\ Math.\ 318, 156-177 (1980)


\bibitem{vdP84}M.~van der Put. {\em Stable Reductions of Algebraic Curves}, Indag.\ Math.\ 46, 461--478 (1984) 

\bibitem{vdPV01}M.~van der Put, H.~Voskuil. {\em Discontinuous subgroups of $\PGL_2(K)$}, preprint 
(2001)


\end{thebibliography}
%

\subsubsection*{Acknowledgements}
The author thanks F.~Kato for distributing his preprint \cite{Kat00},
and Y.~Andr\'e for mentioning the existence of \cite{vdPV01}
after reading a first draft of \cite{Brad02}.
He is also indebted to his thesis advisor F.~Herrlich for answering 
numerous questions.


%
%


\end{document}